# Ακέραια ισοπεριμετρικά τρίγωνα: Συγκριτική θεώρηση δύο διαφορετικών προσεγγίσεων στην επίλυση του ίδιου προβλήματος

**Τάσος Πατρώνης**
Αφυπηρετήσας Επ. Καθηγητής Πανεπιστημίου Πατρών

**Ιωάννης Ρίζος**
ioarizos@uth.gr
Πανεπιστήμιο Θεσσαλίας, Τμήμα Μαθηματικών

**Περίληψη**

Στο άρθρο αυτό αρχικά εισάγουμε τον αναγνώστη στο ισοπεριμετρικό πρόβλημα για την περίπτωση των ακέραιων τριγώνων, στην ακολουθία του Alcuin και στο πώς αυτή σχετίζεται με το πλήθος των διαφορετικών ακέραιων τριγώνων με δοσμένη περίμετρο. Κατόπιν παρουσιάζουμε και συζητάμε δύο διαφορετικές προσεγγίσεις στο παραπάνω πρόβλημα. Η μία προσέγγιση οφείλεται σε έναν επαγγελματία μαθηματικό, καθηγητή Πανεπιστημίου και χρησιμοποιεί μια γνωστή τεχνική υπολογισμού πληθαρίθμων πεπερασμένων συνόλων, ενώ η άλλη σε έναν μαθητή Λυκείου με ιδιαίτερο ενδιαφέρον για τα Μαθηματικά και είναι καθαρά αριθμοθεωρητική. Κλείνοντας, καταλήγουμε σε κάποια ενδεικτικά συμπεράσματα και συζητάμε ορισμένες σχετικές προοπτικές.

**Λέξεις-κλειδιά:** ακέραια τρίγωνα, ισοπεριμετρικό πρόβλημα, ακολουθία του Alcuin, στρατηγικές επίλυσης προβλήματος

**Abstract**

In this paper we first study the isoperimetric problem in the case of integer triangles, as well as Alcuin's sequence and how it relates to the number of different integer triangles with a given perimeter. We then present and compare two different approaches to the above problem. The first approach is due to a university professor (a working mathematician), who uses a well-known technique of calculating the number of elements of finite sets, while the second one is due to a high school student with a special interest in mathematics and is purely number-theoretic. Furthermore, some related instructor perspectives are discussed.

### Εισαγωγή

Ένα τρίγωνο λέγεται "ακέραιο" όταν τα μήκη των πλευρών του είναι όλα ακέραιοι αριθμοί. Μια τριάδα φυσικών αριθμών μπορεί να αντιστοιχεί στα μήκη των πλευρών ενός ακέραιου τριγώνου, εφόσον ικανοποιεί την τριγωνική ανισότητα. Άρα ο αριθμός των διαφορετικών ακέραιων τριγώνων με δοσμένη περίμετρο $p$ είναι το πλήθος των διαμερίσεων του $p$ σε τρία θετικά μέρη που ικανοποιούν την τριγωνική ανισότητα. Αυτός ο αριθμός, ας τον συμβολίσουμε με $T(p)$, αποδεικνύεται ότι είναι:

$$T(p) = \begin{cases} \left[\dfrac{p^2}{48}\right], & \text{αν } p \text{ άρτιος} \\ \left[\dfrac{(p+3)^2}{48}\right], & \text{αν } p \text{ περιττός} \end{cases}$$

όπου $[x]$ είναι η *πλησιέστερη ακέραια συνάρτηση* (nearest integer function) στον αριθμό $x$ και η γραφικής της παράσταση φαίνεται στο ακόλουθο σχήμα (Andrews, 1979; Honsberger, 1985; Jenkyns & Muller, 2000; Jordan et al. 1979):

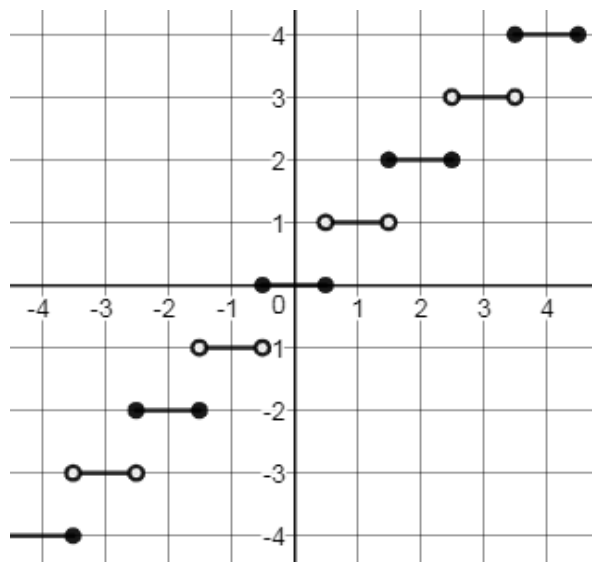

Η ακολουθία $T(p)$ του αριθμού των ακέραιων τριγώνων με περίμετρο $p$, ξεκινώντας από $p = 0$, είναι: 0, 0, 0, 1, 0, 1, 1, 2, 1, 3, 2, 4, 3, 5, 4, 7, 5, 8, 7, 10, 8, 12, 10, 14, 12,... Η παραπάνω ακολουθία είναι γνωστή ως *ακολουθία του Alcuin* (Bindner & Erickson, 2012) και οι όροι της συμπίπτουν με τους συντελεστές του $x^p$ στο ανάπτυγμα της δυναμοσειράς

$$\sum_{p=0}^{\infty} T(p)x^p = \frac{x^3}{(1-x^2)(1-x^3)(1-x^4)}.$$

Η ακολουθία του Alcuin, προκύπτει από τους συντελεστές του αναπτύγματος της παρακάτω συνάρτησης σε μια τυπική (formal) σειρά Maclaurin:

$$\frac{x^3}{(1-x^2)(1-x^3)(1-x^4)} = 1 + x^2 + x^3 + 2x^4 + x^5 + \cdots$$

Γιατί όμως το πλήθος όλων των δυνατών τριγώνων με ακέραιες πλευρές και δοσμένη περίμετρο $p$ ισούται με τον συντελεστή του $x^p$ στο ανάπτυγμα της δυναμοσειράς

$$(A) \quad \frac{x^3}{(1-x^2)(1-x^3)(1-x^4)} \, ;$$

Η απάντηση είναι η εξής: Εφόσον ισχύει

$$\frac{x}{1-x^2} = x + x^3 + x^5 + x^7 \dots$$
$$\frac{x}{1-x^3} = x + x^4 + x^7 + x^{10} \dots$$
$$\frac{x}{1-x^4} = x + x^5 + x^9 + x^{13} \dots,$$

το τυπικό γινόμενο των τριών παραπάνω ισοτήτων είναι το άθροισμα των $a(p)x^p$ για όλα τα $p \geq 3$, όπου $a(p)$ είναι το πλήθος των τρόπων με τους οποίους το $p$ μπορεί να παρασταθεί ως άθροισμα ενός αριθμού της μορφής $2k+1$, ενός αριθμού της μορφής $3k+1$ και ενός αριθμού της μορφής $4k+1$ (με $k$ ακέραιο $\geq 1$). Αυτό είναι ίσο με το πλήθος των τρόπων με τους οποίους το $p-3$ μπορεί να αναπαρασταθεί ως άθροισμα ενός μη αρνητικού πολλαπλάσιου του 2, ενός μη αρνητικού πολλαπλάσιου του 3 και ενός μη αρνητικού πολλαπλάσιου του 4.

Επίσης, για να είναι τρεις θετικοί ακέραιοι αριθμοί $α, β, γ$ πλευρές ενός τριγώνου, πρέπει να ισχύει η τριγωνική ανισότητα: $α < β + γ$, $β < γ + α$, $γ < α + β$. Υποθέτοντας ότι $1 \leq α \leq β \leq γ$, οι δύο πρώτες ανισότητες ικανοποιούνται αυτόματα, όμως πρέπει να έχουμε $γ < α + β - 1$. Θέτοντας $x = β - α$ και $y = γ - β$, γράφουμε την τελευταία ανισότητα ως $α + x + y \leq α + α + x - 1$, δηλαδή $y + 1 \leq α$, ενώ η περίμετρος του τριγώνου είναι $α + β + γ = 3α + 2x + y$. Έτσι το πλήθος των τριγώνων με περίμετρο $p$ είναι το πλήθος των διατεταγμένων ακέραιων τριάδων $(α, x, y)$ με $x \geq 0$, $y \geq 0$, $y + 1 \leq α$ και $3α + 2x + y = p$. Ισοδύναμα, θέτοντας $z = α - y - 1$, το πλήθος των διατεταγμένων τριάδων $(x, y, z)$ μη αρνητικών ακεραίων με $2x + 3z + 4y = p - 3$. Όμως αυτός, όπως είδαμε παραπάνω, είναι ακριβώς ο συντελεστής του $x^p$ στο ανάπτυγμα της δυναμοσειράς $(A)$.

### Περιγραφή της δικής μας δράσης

Καθώς το ισοπεριμετρικό πρόβλημα, ακόμα και στην περίπτωση των ακέραιων τριγώνων, είναι διάσημο (Blåsjö, 2005; East & Niles, 2019; Krier & Manvel, 1998), θελήσαμε να διαπιστώσουμε με ποιον τρόπο θα το αντιμετώπιζε εκ νέου η μαθηματική και εκπαιδευτική κοινότητα στην Ελλάδα. Για τον λόγο αυτό, το καλοκαίρι του 2022 στείλαμε μέσω ηλεκτρονικού ταχυδρομείου πρόσκληση σε: α) καθηγητές Μαθηματικών Δευτεροβάθμιας και Τριτοβάθμιας Εκπαίδευσης, και β) μαθητές οι οποίοι είχαν πάρει μέρος σε διαγωνισμούς της Ε.Μ.Ε. ή γενικότερα επιδείκνυαν έμπρακτο ενδιαφέρον για τα Μαθηματικά. Καθώς υπήρχε μεγάλη καθυστέρηση στις απαντήσεις, απευθυνθήκαμε ιδιαίτερα στον καθηγητή του Τμήματος Μαθηματικών του Πανεπιστημίου Κρήτης Νικόλαο Τζανάκη και τον μαθητή Αστέριο Βαρσάμη-Κυρατλίδη από τη Θεσσαλονίκη, ο οποίος μόλις είχε τελειώσει τη Β΄ Λυκείου. Και οι δύο ανταποκρίθηκαν σύντομα, με πλήρεις λύσεις του προβλήματος, σε ένα τουλάχιστον από τα ερωτήματα.

Το άρθρο αυτό δεν είναι μια έρευνα καταγραφής και ανάλυσης μαθηματικών αντιλήψεων ή αυθόρμητων ευρετικών επίλυσης προβλημάτων. Μεθοδολογικά, μοιάζει αρκετά με τις επώνυμες συνεντεύξεις στον Τύπο ή με την πρόσκληση του Hadamard (1945) σε γνωστούς επιστήμονες του καιρού του να εκθέσουν τον τρόπο σκέψης τους στο βιβλίο *The Psychology of Invention in the Mathematical Field*. Στην πραγματικότητα οι συμμετέχοντες στην έρευνα είναι, εμμέσως, συν-συγγραφείς του άρθρου μας, γι' αυτό δεν έχει νόημα να μην αποκαλύψουμε τα ονόματά τους.

### Το πρόβλημα

Το πρόβλημα διατυπώθηκε ως εξής:

*Έστω ότι $α, β, γ$ και $p$ είναι φυσικοί αριθμοί $\geq 1$, με $α + β + γ = p$.*
*Ζητείται:*
  **α.** *Αν $α, β, γ$ πλευρές τριγώνου, να βρεθεί το πλήθος των μη διατεταγμένων τριάδων $(α, β, γ)$ συναρτήσει του $p$.*
  **β.** *Αν $Ε$ το εμβαδόν τριγώνου με πλευρές μήκους $α, β, γ$, να βρεθεί η μέγιστη τιμή του $Ε$ συναρτήσει του $p$.*

### Η απάντηση του μαθητή στα ερωτήματα α και β

**α.** Θα υπολογίσω το πλήθος των μη διατεταγμένων τριάδων $(α, β, γ)$ για τις διάφορες μορφές του $p \bmod 12$ για να αποφύγω τη χρήση του ακέραιου μέρους.

Αφού οι τριάδες δεν είναι διατεταγμένες, θεωρώ ότι $α \geq β \geq γ$.
Από την τριγωνική ανισότητα προκύπτει ότι $2α < α + β + γ = p$.
Άρα, $α \leq \frac{p-2}{2}$ αν $p$ άρτιος και $α \leq \frac{p-1}{2}$ αν $p$ περιττός.
Επιπλέον ισχύει, $p = α + β + γ \leq 3α$.
Άρα, $α \geq \frac{p}{3}$ αν $p \equiv 0 \bmod 3$, $α \geq \frac{p+2}{3}$ αν $p \equiv 1 \bmod 3$ και $α \geq \frac{p+1}{3}$ αν $p \equiv 2 \bmod 3$.

Βρίσκω όλες τις δυνατές τριάδες με κατασκευαστικό τρόπο.
Για $p = 12n + 1$, όλες οι τριάδες είναι:
$(6n, 6n, 1)$,     $(6n, 6n − 1, 2)$,    …, $(6n, 3n + 1, 3n)$.
(πλήθος τριάδων: $3n$)
$(6n − 1, 6n − 1, 3)$,  $(6n − 1, 6n − 2, 4)$,  …, $(6n − 1, 3n + 1, 3n + 1)$.
(πλήθος τριάδων: $3n − 1$)
$(6n − 2, 6n − 2, 5)$,  $(6n − 2, 6n − 3, 6)$,  …, $(6n − 2, 3n + 2, 3n + 1)$.
(πλήθος τριάδων: $3n − 3$)
………………………………
$(4n + 1, 4n + 1, 4n − 1), (4n + 1, 4n, 4n)$
(πλήθος τριάδων: 2).

Έστω το σύνολο
$A = \{(6n, 3n + 1, 3n), (6n − 1, 3n + 1, 3n + 1), … , (4n + 1, 4n, 4n)\}$.
Αν η τριάδα $(α, β, β)$ με $α ≥ β + 2$ ανήκει στο σύνολο $Α$, τότε και η τριάδα $(α − 1, β + 1, β)$ ανήκει στο σύνολο $Α$, αφού αν ισχύει η τριγωνική ανισότητα για την τριάδα $(α, β, β)$, τότε (αποδεικνύεται εύκολα ότι) θα ισχύει και για την τριάδα $(α − 1, β + 1, β)$. Η επόμενη τριάδα που θα ανήκει και αυτή στο σύνολο $Α$ είναι η $(α − 2, β − 1, β − 1)$. Άρα, ο μικρότερος όρος της τελευταίας τριάδας κάθε σειράς διαφέρει κατά 1 ανά δύο σειρές.

Είναι προφανές ότι ο μικρότερος όρος της πρώτης τριάδας κάθε σειράς διαφέρει κατά 2 ανά μία σειρά. Είναι ακόμα προφανές ότι το πλήθος των τριάδων κάθε σειράς υπολογίζεται από τη διαφορά του μικρότερου όρου της τελευταίας τριάδας κάθε σειράς από τον μικρότερο όρο της πρώτης τριάδας κάθε σειράς, προσθέτοντας, σε αυτή τη διαφορά μια μονάδα. Από αυτά προκύπτει ότι το πλήθος των τριάδων είναι το άθροισμα όλων των αριθμών από το 1 μέχρι το $3n$ αφαιρώντας κάθε αριθμό $x$ με $1 ≤ x ≤ 3n$, όπου $x ≡ 1 mod 3$. Έτσι, για $p = 12n + 1$ το πλήθος των τριάδων είναι $3n^2 + 2n$.[1]

Με όμοιο τρόπο προκύπτουν και οι υπόλοιποι τύποι, που μας δίνουν την πλήρη απάντηση στο πρώτο ερώτημα:
Για $p = 12n$, το πλήθος των τριάδων είναι $3n^2$
Για $p = 12n + 1$, το πλήθος των τριάδων είναι $3n^2 + 2n$
Για $p = 12n + 2$, το πλήθος των τριάδων είναι $3n^2 + n$
Για $p = 12n + 3$, το πλήθος των τριάδων είναι $3n^2 + 3n + 1$

---

[1] Ο μαθητής εννοεί εδώ ότι, από το άθροισμα $S_{3n} = \frac{3n(3n+1)}{2}$ θα πρέπει να εξαιρεθούν οι όροι $1, 4, 7, …, 3n − 2$, που είναι $n$ το πλήθος και έχουν άθροισμα $S_n = \frac{n(3n−1)}{2}$.

Για $p = 12n + 4$, το πλήθος των τριάδων είναι $3n^2 + 2n$
Για $p = 12n + 5$, το πλήθος των τριάδων είναι $3n^2 + 4n + 1$
Για $p = 12n + 6$, το πλήθος των τριάδων είναι $3n^2 + 3n + 1$
Για $p = 12n + 7$, το πλήθος των τριάδων είναι $3n^2 + 5n + 2$
Για $p = 12n + 8$, το πλήθος των τριάδων είναι $3n^2 + 4n + 1$
Για $p = 12n + 9$, το πλήθος των τριάδων είναι $3n^2 + 6n + 3$
Για $p = 12n + 10$, το πλήθος των τριάδων είναι $3n^2 + 5n + 2$
Για $p = 12n + 11$, το πλήθος των τριάδων είναι $3n^2 + 7n + 4$

**β.** Για λόγους συντομίας, συμβολίζω το μέγιστο εμβαδόν ως $maxE$, το τρίγωνο με πλευρές μήκους $α, β, γ$ ως $(α, β, γ)$, και αποκαλώ "εύρος" της τριάδας $(α, β, γ)$ την απόλυτη τιμή της διαφοράς του μέγιστου και του ελάχιστου στοιχείου της τριάδας. Προφανώς, το εύρος της τριάδας $(α, α, α)$ θα είναι ίσο με 0.

*Θα δείξω ότι η τριάδα $(α', β', γ')$ δίνει μεγαλύτερο εμβαδόν από την τριάδα $(α, β, γ)$ αν και μόνο αν η τριάδα $(α', β', γ')$ έχει μικρότερο εύρος από την τριάδα $(α, β, γ)$.*

Έστω $α \geq β \geq γ$ και $α' \geq β' \geq γ'$, και χωρίς βλάβη της γενικότητας αρκεί να το δείξω για $β = β'$, καθώς έτσι μελετώ το εμβαδόν ενός τριγώνου με σταθερή τη βάση και την περίμετρο, μεταβάλλοντας τα μήκη των δύο άλλων πλευρών του.

*Απόδειξη*: Από τον τύπο του Ήρωνα για το εμβαδόν του τριγώνου, έχουμε:
$E = \sqrt{τ(τ-α)(τ-β)(τ-γ)}$, με $τ = p/2$ και
$E' = \sqrt{τ(τ-α')(τ-β')(τ-γ')}$, με $τ = p/2$.
Αν $E < E'$, δηλ. $\sqrt{τ(τ-α)(τ-β)(τ-γ)} < \sqrt{τ(τ-α')(τ-β')(τ-γ')}$, τότε αφού $β = β'$, προκύπτει: $(τ-α)(τ-γ) < (τ-α')(τ-γ')$, άρα $αγ < α'γ'$, καθώς $α + γ = p - β = p - β' = α' + γ'$.

Ακόμα, αφού $α + γ = α' + γ' \Leftrightarrow α^2 + 2αγ + γ^2 = α'^2 + 2α'γ' + γ'^2$.

Ακόμα $αγ < α'γ'$, άρα $-4αγ > -4α'γ'$. Από τα δύο τελευταία με πρόσθεση κατά μέλη προκύπτει: $(α - γ)^2 > (α' - γ')^2 \Leftrightarrow α - γ > α' - γ'$. Άρα, το συμπέρασμά μας ισχύει. Με όμοιο τρόπο αποδεικνύεται και το αντίστροφο. Άρα, η τριάδα με το ελάχιστο δυνατό εύρος, μεγιστοποιεί το εμβαδόν του τριγώνου.

Είναι προφανές ότι θα υπάρχει τριάδα με εύρος ίσο με 0 αν και μόνο αν $p \equiv 0 \bmod 3$. Αν $p = 3n$, υπάρχει η τριάδα $(α', β', γ') = (n, n, n)$ με εύρος ίσο με 0. Άρα, για $p \equiv 0 \bmod 3$, είναι $maxE = \frac{p^2\sqrt{3}}{36}$.

Αν $p = 3n + 1$, υπάρχει τριάδα $(α', β', γ') = (n + 1, n, n)$ με εύρος ίσο με 1, που είναι το ελάχιστο δυνατό. Άρα, για $p \equiv 1 \bmod 3$, είναι $maxE = \frac{p+2}{12} \cdot \sqrt{\frac{p-4}{3}} \cdot \sqrt{p}$.

Αν $p = 3n + 2$, υπάρχει τριάδα $(α', β', γ') = (n + 1, n + 1, n)$ με εύρος ίσο με 1, που είναι το ελάχιστο δυνατό. Άρα, για $p \equiv 2 \bmod 3$, είναι $maxE = \frac{p-2}{12} \cdot \sqrt{\frac{p+4}{3}} \cdot \sqrt{p}$.

Οι τρεις αυτοί τύποι που απαντούν στο δεύτερο ερώτημα, συνοψίζονται στον τύπο:

$$maxE = \frac{p + 2v}{12} \sqrt{\frac{p - 4v}{3}} \sqrt{p}, με\ p \equiv v \bmod 3\ και\ v \in \{-1, 0, 1\}.$$

**Συνέντευξη με τον μαθητή σχετικά με τον τρόπο σκέψης του**

Θέλοντας να διερευνήσουμε σε μεγαλύτερο βάθος την πορεία σκέψης του μαθητή κατά την επίλυση του προβλήματος, πήραμε από αυτόν την ακόλουθη συνέντευξη (οι υπογραμμίσεις είναι δικές μας):

**Ερώτηση:** Πώς ακριβώς σκέφθηκες τη συγκεκριμένη λύση; Αναφέρομαι στο α) ερώτημα, δηλαδή στην εύρεση του πλήθους των μη διατεταγμένων τριάδων α, β, γ, συναρτήσει της περιμέτρου p. Γιατί επέλεξες το $p \bmod 12$, έναντι π.χ. μιας λύσης που θα ενέπλεκε το ακέραιο μέρος ή την πλησιέστερη ακέραια συνάρτηση;

**Απάντηση:** Στην προσπάθεια να φτιάξω μια κατασκευή, αυτή του συνόλου Α, ώστε να μπορώ να καταμετρήσω εύκολα το πλήθος, των τρόπων, είχα ξεκινήσει ερευνώντας για $p \bmod 2$. Αυτό γιατί έκανα την παρατήρηση ότι αν $α_n$ το πλήθος των τρόπων για $p = n$, τότε $a_{2n+1} = a_{2n+4}$. Την παρατήρηση την έκανα έχοντας μετρήσει με καταγραφή το πλήθος των τρόπων για κάθε $p \leq 22$. Σταμάτησα στο 22 τυχαία, ελπίζοντας ότι θα είναι αρκετό. Την κατασκευή που υπάρχει στη λύση ξεκίνησα να την φτιάχνω για $p = 2n + 1$, όμως στο τέλος της πρώτης κιόλας σειράς εμφανίστηκε το ακέραιο μέρος του $n/2$. Οπότε προχώρησα για $p \bmod 4$.

Συνεχίζοντας την κατασκευή για $p = 4n + 1$, στην τελευταία τριάδα της τελευταίας σειράς εμφανίστηκε το ακέραιο μέρος του $(4n + 1)/3$. Έτσι κατέληξα στο να επιλέξω $p \bmod 12$, το οποίο, τελικά, ήταν αρκετό για να γίνει η κατασκευή χωρίς κανένα ακέραιο μέρος. Ο λόγος που ήθελα να μην χρησιμοποιήσω το ακέραιο μέρος στην κατασκευή ήταν το ότι ήθελα να είμαι βέβαιος ότι δεν θα κάνω κάποιο λάθος, όντας, τότε, όχι αρκετά εξοικειωμένος με αυτό. *Ειδικά την πλησιέστερη ακέραια συνάρτηση, που φέρατε ως παράδειγμα, δεν την είχα συναντήσει ή χρησιμοποιήσει ποτέ.*

**Ερώτηση:** Στο δεύτερο ερώτημα, δηλαδή στην εύρεση του μέγιστου εμβαδού συναρτήσει της περιμέτρου $p$, παρατήρησες αν υπάρχει κάποια γεωμετρική ερμηνεία ή σκέφθηκες καθαρά αριθμοθεωρητικά;

**Απάντηση:** *Πλέον, καταλαβαίνω* ότι η γεωμετρική ερμηνεία αυτού που γράφω είναι ότι ελέγχω πότε μεγιστοποιείται το εμβαδόν τριγώνου με σταθερή βάση και η τρίτη κορυφή του να κινείται πάνω σε έλλειψη με εστίες τις δύο σταθερές κορυφές, το οποίο συμβαίνει όταν η τρίτη κορυφή βρίσκεται όσο το δυνατόν κοντινότερα στο σημείο τομής της έλλειψης με τη μεσοκάθετο της βάσης του τριγώνου. *Τότε όμως, το πρώτο που σκέφτηκα*, για να χρησιμοποιήσω τα μήκη πλευρών τριγώνου με γνωστή περίμετρο, ήταν ο τύπος του Ήρωνα. Στη συνέχεια, με δοκιμές σε μικρούς, εύκολα υπολογίσιμους, αριθμούς για τα μήκη των πλευρών του τριγώνου, κατέληξα στην παρατήρηση που έγραψα, απέδειξα και χρησιμοποίησα στη λύση μου. Οπότε, αξιοποιώντας αυτό το λήμμα, δηλαδή, ότι για το ελάχιστο δυνατό εύρος της τριάδας (α, β, γ) το εμβαδόν του τριγώνου πλευρών α, β, γ μεγιστοποιείται και με εύρος της τριάδας (α, β, γ) εννοούμε την απόλυτη τιμή της διαφοράς του μέγιστου και του ελάχιστου στοιχείου της τριάδας, η απάντηση προέκυψε αμέσως.

### Η απάντηση του καθηγητή στο ερώτημα α

Έστω δοσμένος θετικός ακέραιος $p$ (περίμετρος τριγώνου). Ζητούμε το πλήθος των τριάδων $(a, b, c) \in \mathbb{N}^3$ που ικανοποιούν τις συνθήκες
$$a \leq b \leq c < a + b, \quad a + b + c = p. \quad (1)$$
Προφανώς, $c \geq p/3$. Επίσης $p - c = a + b > c$, άρα $c < p/2$. Έτσι,
$$\frac{p}{3} \leq c < \frac{p}{2}. \quad (2)$$
Για κάθε σταθερό $c$ στο παραπάνω διάστημα, οι σχέσεις $a + b = p - c$ και $a \leq b \leq c$ δίνουν
$$a = p - c - b \geq p - 2c, \quad a = p - c - b \leq p - c - a, \quad (3)$$
οπότε

$$p - 2c \le a \le \left\lfloor \frac{p-c}{2} \right\rfloor.$$

Άρα το πλήθος $T$ των $(a, b, c)$ που ικανοποιούν τις (1) είναι, το πολύ, όσο το πλήθος των $(c, a)$ που ικανοποιούν τις (2) και (3). Αλλά και αντίστροφα, σε κάθε τέτοιο ζευγάρι $(a, c)$ αντιστοιχεί μια τριάδα $(a, b, c)$ που ικανοποιεί τις (1). Πράγματι, αν τα $c, a$ ικανοποιούν τις (2) και (3) και θέσουμε $b = p - (a + c)$, τότε $p - c - b = a \ge p - 2c$ άρα $c \ge b$. Επίσης, $p - c - b = a \le (p - c)/2$, άρα $2p - 2c - 2b \le p - c$, οπότε $p - c \le 2b$. Όμως $p - c \ge 2a$, άρα $2b \ge 2a$. Έτσι, $b \ge a$. Τέλος, $a + b = p - c > c$ λόγω της (2), άρα $a + b > c$ και, εξ ορισμού του $b$, $a + b + c = p$.

Συνεπώς, το $T$ *ισούται* με το πλήθος των $(c, a)$ που ικανοποιούν τις συνθήκες (2) και (3), άρα

$$T = \sum_{c=\left\lceil\frac{p}{3}\right\rceil}^{\left\lfloor\frac{p-1}{2}\right\rfloor} \sum_{a=p-2c}^{\left\lfloor\frac{p-c}{2}\right\rfloor} 1 = \sum_{c=\left\lceil\frac{p}{3}\right\rceil}^{\left\lfloor\frac{p-1}{2}\right\rfloor} \left( \left\lfloor \frac{p-c}{2} \right\rfloor - (p - 2c) + 1 \right).$$

### Συζήτηση – Κάποια ενδεικτικά συμπεράσματα

Οι λύσεις του παραπάνω προβλήματος προέρχονται από εκπροσώπους δύο διαφορετικών ομάδων λυτών: από το ένα μέρος πρόκειται για τους επαγγελματίες μαθηματικούς και ερευνητές, και από το άλλο για τους μαθητές με έμπρακτο ενδιαφέρον για τα Μαθηματικά. Οι δύο αυτοί εκπρόσωποι φαίνεται πως εφαρμόζουν, τουλάχιστον στην έρευνά μας, δύο διαφορετικές στρατηγικές λύσεων. Η πρώτη προσέγγιση χρησιμοποιεί τη γνωστή τεχνική (Μαμωνά-Downs & Παπαδόπουλος, 2017) της κατασκευής αντιστοιχίας "ένα-προς-ένα" για τον υπολογισμό των πληθάριθμων πεπερασμένων συνόλων. Η δεύτερη προσέγγιση είναι καθαρά αριθμοθεωρητική. Όμως, ίσως από αφοσίωση στην Αριθμοθεωρία, ο μαθητής δεν συνειδητοποιεί μια γεωμετρική εικόνα που είναι συμφυής με τη λύση του: μια έλλειψη με εστιακή απόσταση ίση με την πλευρά $\beta$ του τριγώνου και σημεία, τέτοια ώστε το άθροισμα των αποστάσεών τους από τις δύο εστίες να είναι ίσο με $α + γ = p - β = σταθ.$, όπου $α, γ$ οι άλλες δύο πλευρές του τριγώνου.

Θα θέλαμε να επισημάνουμε πως η λύση του μαθητή υποδεικνύει μια αισιόδοξη προοπτική· το τι μπορούν να πετύχουν οι μαθητές που αγαπούν τα Μαθηματικά όταν τους ενθαρρύνουμε να λύσουν "ασυνήθιστα" προβλήματα. Ταυτόχρονα όμως δείχνει και το τι χάνουμε αφαιρώντας από την ύλη του Λυκείου τη Θεωρία Αριθμών και υποβαθμίζοντας το μάθημα της Γεωμετρίας, πόσο πιο φτωχή γίνεται η μαθηματική μας εκπαίδευση, καθώς και η εκπαίδευση των νέων ερευνητών στη διεξαγωγή πρωτότυπης έρευνας σε βασικούς κλάδους των Μαθηματικών. Ιδιαίτερα για τους λεγόμενους "χαρισματικούς" μαθητές στα Μαθηματικά, ίσως η καλύτερη προοπτική δεν είναι η ενασχόληση με τα ειδικής φύσης προβλήματα των μαθηματικών διαγωνισμών, αλλά μια ευρύτερη παιδεία που θα περιλαμβάνει μεταξύ άλλων τη βασική Θεωρία Αριθμών καθώς και πρωτότυπα προβλήματα που συνδυάζουν τη Γεωμετρία με την Άλγεβρα.

Μια μελλοντική έρευνα, σε επίπεδο Τριτοβάθμιας ή/και Δευτεροβάθμιας Εκπαίδευσης, θα μπορούσε να μελετήσει σε βάθος τις στάσεις και τις στρατηγικές των φοιτητών/μαθητών με ιδιαίτερο ενδιαφέρον για τα Μαθηματικά πάνω σε παρόμοια προβλήματα και να κάνει προτάσεις για τη θέση, το περιεχόμενο και τον ρόλο της Θεωρίας Αριθμών στα Προγράμματα Σπουδών στο Πανεπιστήμιο και το Λύκειο.

**Αναφορές**
Andrews, G. E. (1979). A Note on Partitions and Triangles with Integer Sides. *The American Mathematical Monthly*, *86*(6), 477-478.
Bindner, D. J. & Erickson, M. (2012). Alcuin's Sequence. *The American Mathematical Monthly*, *119*(2), 115-121.
Blåsjö, V. (2005). The Isoperimetric Problem. *The American Mathematical Monthly*, *112*(6), 526-566.
East, J., & Niles, R. (2019). Integer Triangles of Given Perimeter: A New Approach via Group Theory. *The American Mathematical Monthly*, *126*(8), 735-739.
Hadamard, J. (1945). *The psychology of invention in the mathematical field*. Princeton: Princeton University Press.
Honsberger, R. (1985). *Mathematical Gems III*. Washington, DC: The Mathematical Association of America.
Iverson, K. E. (1962). *A Programming Language*. New York: Wiley.
Jenkyns, T., & Muller, E. (2000). Triangular Triples from Ceilings to Floors. *The American Mathematical Monthly*, *107*(7), 634-639.
Jordan, J. H., Walch, R., & Wisner, R. J. (1979). Triangles with Integer Sides. *The American Mathematical Monthly*, *86*(8), 686-689.
Krier, N., & Manvel, B. (1998). Counting Integer Triangles. *Mathematics Magazine*, *71*(4), 291-295.
Μαμωνά-Downs, Γ., & Παπαδόπουλος, Ι. (2017). *Επίλυση προβλήματος στα μαθηματικά*. Ηράκλειο: Πανεπιστημιακές Εκδόσεις Κρήτης.

## ΠΑΡΑΡΤΗΜΑ

Με $\lfloor x \rfloor$ συμβολίζουμε τη *συνάρτηση κατώτατου ορίου* (*floor function*), η οποία δίνει τον μεγαλύτερο ακέραιο αριθμό που είναι μικρότερος ή ίσος του $x$ (βλ. Iverson, 1962, σ. 12).

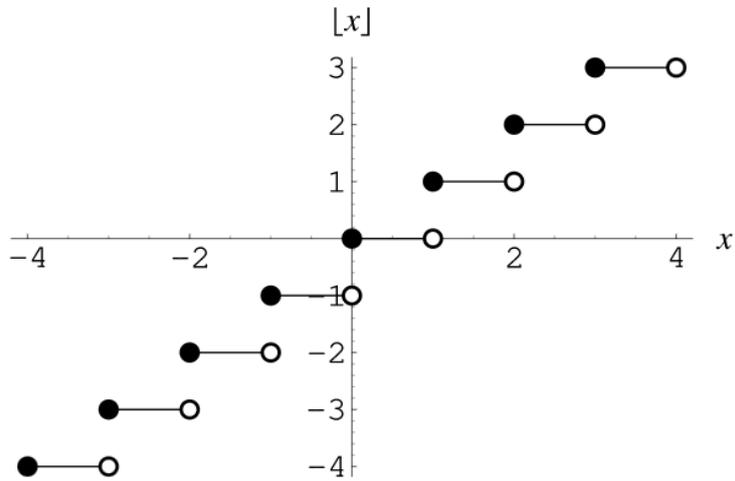

Με $\lceil x \rceil$ συμβολίζουμε τη *συνάρτηση ανώτατου ορίου* (*ceiling function*), η οποία δίνει τον μικρότερο ακέραιο αριθμό που είναι μεγαλύτερος ή ίσος του $x$ (βλ. Iverson, 1962, σ. 12).

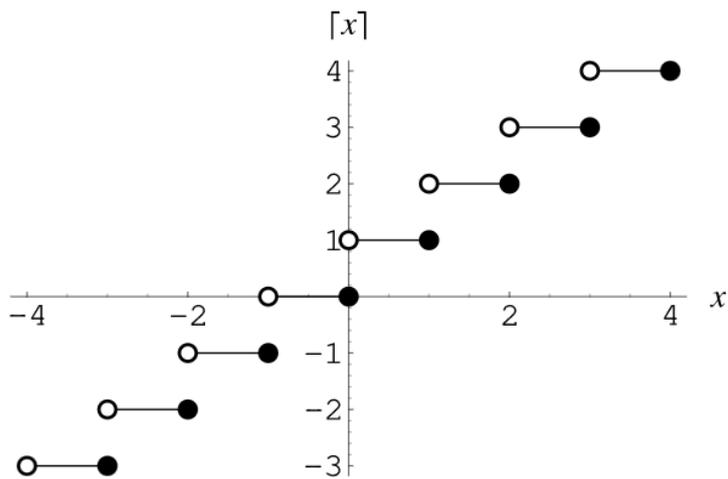